\documentclass[12pt,a4paper,english,intoc,bibliography=totoc,index=totoc,BCOR10mm,captions=tableheading,titlepage,fleqn]{scrbook}
\usepackage{lmodern}

\usepackage[T1]{fontenc}
\usepackage[latin9]{inputenc}
\usepackage{babel}
\usepackage{amsmath}
\usepackage{amssymb}
\usepackage[unicode=true,
 bookmarks=true,bookmarksnumbered=true,bookmarksopen=true,bookmarksopenlevel=2,
 breaklinks=false,pdfborder={0 0 0},backref=false,colorlinks=false]
 {hyperref}
\hypersetup{pdftitle={Abstract},
 pdfauthor={Uwe Stöhr},
 pdfpagelayout=OneColumn, pdfnewwindow=true, pdfstartview=XYZ, plainpages=false}
\usepackage{breakurl}

\makeatletter

\usepackage{calc}

\makeatother

\begin{document}
\addcontentsline{toc}{chapter}{Abstract} 

\noindent \begin{center}
\textbf{\Large A Borel maximal eventually different family}
\par\end{center}{\Large \par}

\noindent \begin{center}
{\large Haim Horowitz and Saharon Shelah}
\par\end{center}{\large \par}

\noindent \begin{center}
\textbf{Abstract}
\par\end{center}

\noindent \begin{center}
We construct a Borel maximal eventually different family.%
\footnote{Date: May 12, 2016.

2000 Mathematics Subject Classification: 03E15, 03E25.

Keywords: Borel, MED families, MAD families.

Publication 1089 of the second author.%
}
\par\end{center}

\textbf{\large Introduction}{\large \par}

Maximal almost disjoint families and their relatives have been studied
by set theorists for decades. As the construction of such families
is typically being done using the axiom of choice, questions about
their definability naturally arise. The definability of mad families
was investigated by Mathias who proved the following:

\textbf{Theorem {[}Ma{]}: }There are no analytic mad families.

The possibility of the non-existence of mad families was investigated
by the authors in {[}HwSh:1090{]} where the following was proved (earlier
such results were proven by Mathias and Toernquist using Mahlo and
inaccessible cardinals, respectively):

\textbf{Theorem {[}HwSh:1090{]}: }$ZF+DC+"$There are no mad families$"$
is equiconsistent with $ZFC$.

In this paper we shall study maximal eventually different families
in $\omega^{\omega}$. Recall that $f,g \in \omega^{\omega}$ are
eventually different if $f(n) \neq g(n)$ for large enough $n$. A
family $\mathcal F \subseteq \omega^{\omega}$ is a maximal eventually
family if the members of $\mathcal F$ are pairwise eventually different,
and $\mathcal F$ is maximal with respect to this property. Our main
goal is to construct in $ZF$ a Borel maximal eventually different
family, thus answering a question asked by several set theorists (see
for example {[}Br{]}, {[}KSZ{]} and {[}To{]}) and showing that the
analog for the above theorems is not true for maximal eventually different
families. While in the current paper we do not attempt to find the
minimum possible Borel complexity of a MED family, we intend to prove
in future work that there exist closed MED families.

In a subsequent paper we shall also prove a similar result for maximal
cofinitary groups and investigate some connections with Borel combinatorics
and large cardinals.

\textbf{\large The proof}{\large \par}

\textbf{Theorem 1 ($ZF$): }There exists a Borel MED family.

Observe that the notion of a Borel MED family can be defined for $A^B$
whenever $|A|=\aleph_0=|B|$, and it's enough to prove that for some
$A$ and $B$ of cardinality $\aleph_0$, there is a Borel MED family
in $A^B$ (with the natural Polish topology).

\textbf{Definition and claim 2: }a. Let $T_*=2^{<\omega}$.

b. $\mathcal{F}_*=\{f : f$ is a function from $2^{<\omega}$ to $H(\aleph_0)\}$.

c. For $n<\omega$ let $\mathcal{F}_n^*=\{f\restriction 2^{<n} : f\in \mathcal{F}_*\}$.

d. For $f,g \in \mathcal{F}_*$ let $eq(f,g)=\{\rho : f(\rho)=g(\rho)\}$
and $dif(f,g)=2^{<\omega} \setminus eq(f,g)$.

e. Let $EDF=\{ \mathcal F \subseteq \mathcal{F}_* :  (\forall f\neq g \in \mathcal F)(|eq(f,g)|<\aleph_0) \}$.

f. Let $MEDF=\{\mathcal F \in EDF : \mathcal F$ is maximal$\}$.

g. Let $\bold B: \mathcal{F}_* \rightarrow 2^{\omega}$ be an injective
continuous function.

h. Let $F_1: \mathcal{F}_* \rightarrow \mathcal{F}_*$ be defined
as $F_1(f)(\rho)=f \restriction 2^{<lg(\rho)}$.

i. Let $G_0=\{F_1(f) : f\in \mathcal{F}_*\}$.

j. Let $G_1$ be the set of $g\in \mathcal{F}_*$ such that for some
$f\in \mathcal{F}_*$, $dif(g,F_1(f))$ is infinite and satisfies:

1. $(\forall n)(|\{ \rho : \bold{B}(f) \restriction  n \nleq \rho \wedge \rho \in dif(g,F_1(f))\}|<\aleph_0)$.

2. For every $\rho \in 2^{<\omega}$, if $\rho \leq \bold{B}(f)$
then there exists at most one $\nu$ such that $\rho \leq \nu \in dif(g,F_1(f))$
and $\nu \cap \bold{B}(f)= \rho$.

k. For $g\in G_1$, let $f_g$ be the unique $f$ as in clause (j).
We shall prove that $f_g$ is indeed unique, and can be Borel-computed
from $f$.

l. For $g\in G_1$ and $f_g$ as above, let $w_g=dif(g,F_1(f_g))$.

m. Let $G_2$ be the set of $g\in G_1$ satisfying (1) and (2) where:

1. $g\restriction w_g=f_g \restriction w_g$.

2. $(\forall \rho \in w_g)(g(\rho) \notin \mathcal{F}_{lg(\rho)}^*)$
or $(\forall \rho \neq \nu \in w_g)(g(\rho) \in \mathcal{F}_{lg(\rho)}^* \wedge g(\rho) \nsubseteq g(\nu))$.

\textbf{Proof (of clause (k)): }Given $g\in G_1$, let $X_1(g)=\{\rho \in T_* : g(\rho) \in \mathcal{F}_{lg(\rho)}^*\}$.
Let $X_2(g)=\{ \rho \in T_* : (\forall \nu_1,\nu_2)(\rho \leq \nu_1 \leq \nu_2 \rightarrow \nu_1,\nu_2 \in X_1(g) \wedge g(\nu_1) \subseteq g(\nu_2))\}$,
$X_3(g)=\{\rho \in  T_* : |\{\nu : \rho \leq \nu \in T_*, \nu \notin X_2(g)\}|<\aleph_0\}$
and $X_4(g)=\{ \rho \in X_3(g) : $ there are no incompatible $\nu_1$
and $\nu_2$ such that $\rho \leq \nu_1,\nu_2 \in T_*$ and $\nu_l \notin X_2(g)$
$(l=1,2)\}$. As $g\in G_1$, there is $f$ as in clause (j). 

We shall now prove that if $\rho \nleq \bold{B}(f)$ then $\rho \in X_3(g)$
and moreover, $\rho \in X_4(g)$: By the definition of $G_1$, $\Lambda_n:=\{ \nu \in T_* : \bold{B}(f) \restriction n \nleq \nu, g(\nu) \neq F_1(f)(\nu)\}$
is finite for every $n<\omega$. Now let $\rho \in T_*$ such that
$\rho \nleq \bold{B}(f)$ and choose a minimal $n$ such that $\bold{B}(f) \restriction n \nleq \rho$.
For every $\rho \leq \nu \in T_*$, if $\nu \notin \Lambda_n$ then
$g(\nu)=F_1(f)(\nu)$, therefore, $\rho \leq \nu_1 \leq \nu_2 \in T_* \wedge \nu_1,\nu_2 \notin \Lambda_n \rightarrow g(\nu_1)=F_1(f)(\nu_1) \subseteq F_1(f)(\nu_2)=g(\nu_2)$.
It follows that $\rho \in X_3(g)$, moreover, by 2(j)(2), $\rho \in X_4(g)$:
There is at most one $\nu$ such that $\rho \leq \nu$ and $\nu \in dif(g,F_1(f))$.
For every $\rho \leq \nu'$ which is not $\leq \nu$, $g(\nu')=f\restriction 2^{<lg(\nu')}$,
hence $\nu' \in X_2(g)$. It follows that $\rho \in X_4(g)$. 

Therefore, for every $n$, $|\{ \rho \in T_* : lg(\rho)=n, \rho \notin X_4(g)\}| \leq 1$.
Note that $X_i(g)$ $(i=1,2,3,4)$ can be simply computed. 

Note that by 2(j)(2), for every $\rho \in 2^{<\omega}$ there exists
$\rho' \in eq(g,F_1(f))$ above it, hence, if $\rho \in X_2(g)$ then
$\rho \in eq(g,F_1(f))$. Now suppose that $\nu_1 \neq \nu_2 \in 2^n \cap dif(g,F_1(f))$.
If $\nu_1 \cap \nu_2 \nleq \bold{B}(f)$, then $\nu_1 \cap \nu_2 \in X_4(g)$,
contradicting the fact that $\nu_1,\nu_2 \notin X_2(g)$ are incomparable.
If $\nu_1 \cap \bold{B}(f)=\nu_2 \cap \bold{B}(f)=\nu_1 \cap \nu_2$,
then we get a contradction to 2(j)(2). The only possibility left is
that $\nu_1 \cap \nu_2 \leq \bold{B}(f)$ but $\nu_1 \cap \bold{B}(f) \neq \nu(2) \cap \bold{B}(f)$,
so wlog $\nu_1 \cap \bold{B}(f)< \nu_2 \cap \bold{B}(f)$. Therefore,
there are at most $n$ elements $\nu \in 2^n$ such that $\nu \in dif(g,F_1(f))$.
As $2^{n-1}>n$ for $3\leq n$, we have established the following:

$(*)$ If $3\leq n$, then for most $\nu \in 2^n$, $g(\nu)=f\restriction 2^{<n}$.

It follows that if $g\in G_1$ then $f_g$ is uniquely determined,
and there exists a Borel function $\bold{B}': \mathcal{F}_* \rightarrow \mathcal{F}_*$
such that $g\in G_1 \rightarrow \bold{B}'(g)=f_g$. $\square$

\textbf{Claim 3: }1. If $g_1,g_2 \in G_2$ and $f_{g_1} \neq f_{g_2}$,
then:

a. $eq(g_1,g_2)$ is finite.

b. $w_{g_1} \cap w_{g_2}$ is finite.

c. $eq(g_2,F_1(f_{g_1}))$ is finite.

2. If $g_1 \in G_2$, $f_0 \in \mathcal{F}_*$ and $f_{g_1} \neq f_0$,
then $eq(g_1,F_1(f_0))$ is finite.

\textbf{Proof: }1. As $\bold{B}$ is injective, $\bold{B}(f_{g_1}) \neq \bold{B}(f_{g_2})$,
therefore $\rho:= \bold{B}(f_{g_1}) \cap \bold{B}(f_{g_2}) \in 2^{<\omega}$
and WLOG $\rho \hat{} (l) \leq \bold{B}(f_{g_l})$. By the definition
of $G_1$, $\{\nu \in w_{g_l} : \rho \hat{} (l) \nleq \nu\}$ is finite
for $l=1,2$, therefore $w_{f_{g_1}} \cap w_{f_{g_2}}$ is finite,
which proves clause (b). Now let $n_*$ be such that $f_{g_1} \restriction 2^{<n_*} \neq f_{g_2} \restriction 2^{<n_*}$.
If $\nu \in 2^{<\omega} \setminus w_{g_1} \setminus w_{g_2} \setminus 2^{\leq n_*}$,
then $g_l(\nu)=F_1(f_{g_l})(\nu)$ $(l=1,2)$ by the definition of
$w_{g_l}$. By the choice of $n_*$ and the definition of $F_1$,
$F_1(f_{g_1})(\nu) \neq F_2(f_{g_2})(\nu)$, so $g_1(\nu) \neq g_2(\nu)$.
Note that $|\{ \nu \in w_{g_2} : g_2(\nu)=F_1(f_{g_1})(\nu) \}| \leq 1$:
By the definition of $G_2$, either $g_2(\nu) \notin \mathcal{F}_{lg(\nu)}^*$
for every $\nu \in w_{g_2}$ (in this case, the above set is empty
by the definition of $F_1$) or $\{g_2(\nu) : \nu \in w_{g_2}\}$
are pairwise incomparable with respect to inclusion, and then as $\{F_1(f_{g_1})(\nu) : \nu \in w_{g_2}\}$
form a chain, the above set has cardinality $\leq 1$. Suppose now
that $\nu \in w_{g_2} \setminus w_{g_1}$, then $g_1(\nu)=F_1(f_{g_1})(\nu)$,
and by the above claim, there is at most one $\nu \in w_{g_2} \setminus w_{g_1}$
such that $g_1(\nu)=g_2(\nu)$. Similarly, there is at most one $\nu \in w_{g_1} \setminus w_{g_2}$
such that $g_1(\nu)=g_2(\nu)$. Therefore, $eq(g_1,g_2)$ is finite,
which proves clause (a). Clause (c) follows from (2).

2. By the definition of $G_2$, either $g_1(\nu) \notin \mathcal{F}_{lg(\nu)}^*$
for every $\nu \in w_{g_1}$ (and therefore $w_{g_1} \cap eq(g_1,F_1(f_0))=\emptyset$),
or $\{g_1(\nu) : \nu \in w_{g_1}\}$ are pairwise incomparable (and
then $|w_{g_1} \cap eq(g_1,F_1(f_0))| \leq 1$). If $\nu \notin w_{g_1}$
is long enough, then $g_1(\nu)=F_1(f_{g_1})(\nu)=f_{g_1} \restriction 2^{<lg(\nu)} \neq f_0 \restriction 2^{<lg(\nu)}=F_1(f_0)(\nu)$.
Together we get the desired conclusion. $\square$

\textbf{Definition 4: }Let $H_3=\{ f\in \mathcal{F}_* : $ there is
$g\in G_2$ such that $f_g=f\}$.

\textbf{Definition 5: }Given a formula $\psi(x)$, we say that the
truth value $TV(\psi(f))$ $(f\in \mathcal{F}_*)$ is Borel-computable
if there exists a Borel function $F: \mathcal{F}_* \rightarrow \{0,1\}$
such that $TV(\psi(f))=true$ iff $F(f)=1$. 

The theorem will follow from the following claim together with claim
8:

\textbf{Claim 6: }There is a Borel function $F_3^*$ such that $Dom(F_3^*)=\mathcal{F}_*$,
$f\in H_3 \iff F_3^*(f) \in G_2$ and $f_{F_3^*(f)}=f$ when $f\in H_3$.
As a consequence, $H_3$ is Borel.

\textbf{Definition 7: }Let $G_4:=\{F_3^*(f) : f\in H_3\} \cup \{F_1(f) : f\in \mathcal{F}_* \setminus H_3\}$.

\textbf{Claim 8: }a. $G_4$ is Borel and $G_4 \subseteq G_0 \cup G_2$
(and $G_2 \subseteq G_1$).

b. $G_4 \in EDF$.

c. $G_4 \in MEDF$.

\textbf{Proof of claim 8: }a. The second part of the claim is obvious.
As for the first part, first observe that $f\in G_4$ iff $TV_1(f)=true$
or $TV_2(f)=true$ where:

1. $TV_1(f)=true$ iff $f\in G_0$ and $F_1^{-1}(f) \notin H_3$ (where
$G_0$ was defined in 2(i)).

2. $TV_2(f)=true$ iff $\bold{B}'(f) \in H_3$ and $f=F_3^*(\bold{B}'(f))$
(where $\bold{B}'$ is the Borel function from claim 2(k), which is
defined in the end of the proof of the claim).

Next observe that $TV_1(f)$ is Borel-computable: It's easy to see
that $G_0$ is closed and $F_1^{-1}$ is continuous on $G_0$. As
$H_3$ is Borel, we're done.

$TV_2(f)$ is Borel-computable as well, as $H_3$ and all of the functions
involved are Borel. It follows that $G_4$ is Borel.

b. Suppose that $g_1 \neq g_2 \in G_4$ as witnessed by $f_{g_1}=f_1$
and $f_{g_2}=f_2$. Clearly, $f_1=f_2$ is impossible, as then, if
$f_1 \in H_3$ then $f_2 \in H_3$, hence $g_1=F_3^*(f_1)=F_3^*(f_2)=g_2$,
and similarly, if $f_1,f_2 \notin H_3$, then $g_1=F_1(f_1)=F_1(f_2)=g_2$.
Therefore, $f_1 \neq f_2$. If $f_1,f_2 \in H_3$ then $g_1,g_2 \in G_2$
and by claim 3(1), $eq(g_1,g_2)$ is finite. If $f_1,f_2 \notin H_3$,
then $g_1=F_1(f_1), g_2=F_1(f_2)$, and by the definition of $F_1$,
$eq(g_1,g_2)$ is finite. If $f_1 \in H_3$ and $f_2 \notin H_3$
or vice versa, then $eq(g_1,g_2)$ is finite by 3(2). 

c. Let $f\in \mathcal{F}_*$, we shall find $g \in G_4$ such that
$eq(f,g)$ is infinite. Denote $\bold{B}(f)$ (from 2(g)) by $\eta_f$.
If $f\in H_3$ then $g=F_3^*(f) \in G_4$ is well-defined. By the
definition of $G_2$ and $F_3^*$, $g\restriction w_g=f \restriction w_g$.
By the definition of $G_2$, $w_g$ is infinite. Therefore, we may
assume that $f\notin H_3$.

\textbf{Case I: For every $n$ there is $\nu$ such that $\eta_f \restriction n \leq \nu \in 2^{<\omega}$
and $f(\nu) \notin \mathcal{F}_{lg(\nu)}^*$.} In this case, choose
the $<_*-$least witness $\nu_n$ for every $n$. There is an infinite
set $A\subseteq \omega$ such that $(lg(\nu_n \cap \eta_f) : n\in A)$
is strictly increasing. Let $g=(f \restriction \{\nu_n : n\in A\}) \cup (F_1(f) \restriction (2^{<\omega} \setminus A))$,
it's straightforward to verify that $g\in G_2$ (by the first possibility
in definition 2(m)(2)) and $f=f_g$, which is a conradiction.

\textbf{Case II: Case I fails, but there are $A \in [\omega]^{\omega}$
and $\bar{\nu}=(\nu_n : n\in A)$ such that $\eta_f \restriction n \leq \nu_n$,
$lg(\nu_n \cap \eta_f)=n$ and $(f(\nu_n) : n \in A)$ are pairwise
incomparable. }In this case, we shall derive a contradiction as in
the previous case (using the second possibility in definition 2(m)(2)).
Note that if $n$ exemplifies the failure of case I, then as $(f(\nu_m) : n\leq m \in A)$
are pairwise incomparable, there is at most one $n\leq m \in A$ such
that $f(\nu_m)=F_1(f)(\nu_m)$. If $n \leq n_*$ and $f(\nu_m) \neq F_1(f)(\nu_m)$
for every $n_* \leq m \in A$, then we define $g$ as in the previous
case, with $\{\nu_m : n_* \leq m \in A\}$ here instead of $\{\nu_n : n\in A\}$
there, and we get a contradiction similarly.

\textbf{Case III: $\neg$Case I$\wedge \neg$Case II. }We shall prove
the following statement:

$(*)$ There are $n_*$, $k^*$ and $f_0,...,f_{k_*} \in \mathcal{F}_*$
such that $\eta_f \restriction n_* \leq \nu \rightarrow f(\nu) \in \{f_0 \restriction 2^{<lg(\nu)},...,f_{k_*-1} \restriction 2^{<lg(\nu)}\}$. 

In order to prove $(*)$, assume that it fails and we shall derive
a contradiction to the assumptions of case III.

Let $n_1$ witness the failure of case I, we choose by induction on
$k$ a triple $(\bar{\eta_k},A_k,f_k)$ such that:

a. $\bar{\eta_k}=(\eta_{k,n} : n_1 \leq n \in A_k)$.

b. $\eta_f \restriction n \leq \eta_{k,n}$ but $\eta_f \restriction (n+1) \nleq \eta_{k,n}$.

c. $f(\eta_{k,n}) \notin \{f_l \restriction 2^{<lg(\eta_{k,n})} : l<k\}$.

d. $A_k \subseteq \omega$ is infinite and $(f(\eta_{k,n}) : n_1 \leq n \in A_k)$
is $\subseteq-$increasing.

e. $f_k \in \mathcal{F}_*$ and $f_k=\underset{n_1 \leq n}{\cup}f(\eta_{k,n})$.

Why can we carry the induction? At stage $k$, let $A_k^1=\{n : n_1 \leq n$
and there is $\eta_{k,n}$ satisfying (b)+(c)$\}$. If $A_k^1$ is
finite, then letting $n_*=max(A_k^1)+1$, $(n_*,k-1,f_0,...,f_{k-1})$
are as required in the above statement $(*)$, contradicting thee
assumption that $(*)$ fails. If $A_k^1$ is infinite, we can choose
for every $n\in A_k^1$ an $\eta_{k,n}$ satisfying (b)+(c) (for example,
by taking the $<_*-$minimal such sequence), by Ramsey's theorem there
is an infinite $A_k \subseteq A_k^1$ such that $(f(\eta_{k,n}) : n\in A_k)$
is either $\subseteq-$increasing, $\subseteq-$decreasing or pairwise
incomparable (note that we don't need any form of the axiom of choice
here, as we can carry the argument in a model of the form $L[X]$).
If the elements of $\{f(\eta_{k,n}) : k\in A_n\}$ are pairwise incomparable,
let $w=\{\eta_{k,n} : n\in A_k\}$ and $g=(f\restriction w) \cup (F_1(f) \restriction (2^{<\omega} \setminus w))$.
It's straightforward to verify that $g\in G_2$ and $f_g=f$ (note
that by the pairwise incomparability of the $f(\eta_{k,n})$s, there
is at most one $\eta_{k,n}$ for which $f(\eta_{k,n})=F_1(f)(\eta_{k,n})$).
Therefore, $f\in H_3$, contradicting our assumption. By the choice
of $n_1$, the sequence $(f(\eta_{k,n}) : n\in A_k)$ can't be $\subseteq-$decreasing,
therefore, it's $\subseteq-$increasing. Let $f_k=\cup \{f(\eta_{k,n}) : n\in A_k\}$,
then $f_k \in \mathcal{F}_*$ and $n\in A_k \rightarrow f(\eta_{k,n})=F_1(f_k)(\eta_{k,n})$,
so we've carried the induction.

We shall now get a contradiction by showing that the assumptions of
case II hold: Note that $k_1 \neq k_2 \rightarrow f_{k_1} \neq f_{k_2}$
(by clauses (c) and (e)). Let $B_0=\omega$, choose $l_0$ such that
$f_0 \restriction 2^{\leq l_0} \neq f_1 \restriction 2^{\leq l_0}$.
Therefore, there are $h_0 \in \{0,1\}$ and an infinite set $B_1 \subseteq \omega \setminus \{0,1\}$
such that $\underset{k\in B_1}{\wedge}f_k \restriction 2^{\leq l_0} \neq f_{h_0} \restriction 2^{\leq l_0}$.
Now choose $i_{1,0} \neq i_{1,1} \in B_1$ and $l_1$ such that $f_{i_{1,0}} \restriction 2^{\leq l_1} \neq f_{i_{1,1}} \restriction 2^{\leq l_1}$.
As before, there are $h_1 \in \{0,1\}$ and an infinite set $B_2 \subseteq B_1 \setminus (i_{1,0}+i_{1,1})$
such that $\underset{k\in B_2}{\wedge}f_k \restriction 2^{\leq l_1} \neq f_{i_{1,h_1}} \restriction 2^{\leq l_1}$.
We continue as above and obtain the sets $B=\{h_0<i_{1,h_1}<i_{2,h_2}<...\}$,
$(B_n : n<\omega)$, $((i_{m,0},i_{m,1}) : m<\omega)$ and $(l_m : m<\omega)$.
For every $k\in B$, if $k=i_{m,h_m}$, choose $n_k \in A_k$ such
that $max\{l_m,n_{k-1}\}<n_k$ and let $\nu_{n_k}=\eta_{k,n_k}$ and
$A=\{n_k : k\in B\}$. It's now easy to verify that $A$ and $(\nu_{n_k} : k\in B)$
satisfy the assumptions of case II, but we shall elaborate: We shall
prove that $f(\nu_{k_1})=f(\eta_{k_1,n_{k_1}})=f_{k_1} \restriction 2^{\leq lg(\eta_{k_1,n_{k_1}})}$
and $f(\nu_{k_2})=f(\eta_{k_2,n_{k_2}})=f_{k_2} \restriction 2^{\leq lg(\eta_{k_2,n_{k_2}})}$
are incomparable for $k_1 \neq k_2 \in B$. Suppose that $k_1=i_{m,h_m}$
and $k_2=i_{j,h_j}$ and wlog $m<j$, then $f_{k_1} \restriction 2^{\leq l_m} \neq f_{k_2} \restriction 2^{\leq l_m}$,
therefore $f_{k_1} \restriction 2^{\leq lg(\eta_{k_1,n_{k_1}})} \neq f_{k_2} \restriction 2^{\leq lg(\eta_{k_1,n_{k_1}})}$
and $f_{k_1} \restriction 2^{\leq lg(\eta_{k_2,n_{k_2}})} \neq f_{k_2} \restriction 2^{\leq lg(\eta_{k_2,n_{k_2}})}$,
and therefore $f(\nu_{k_1})$ and $f(\nu_{k_2})$ are incomparable.
This completes the proof of $(*)$.

Now let $n_*,k_*,f_0,...,f_{k_*-1}$ be as in $(*)$, then for every
$n \geq n_1$, there is $l_n<k_*$ such that the set $Y_n=\{ \rho \in 2^{<\omega} : \eta_f \restriction n \leq \rho$,
$\eta_f(n) \neq \rho(n)$ and $f(\rho)=F_1(f_{l_n})(\rho)\}$ is infinite.
Choose $l_*<k_*$ such that $B=\{n : n_1 \leq n, l_n=l_*\}$ is infinite.

Subcase I: $f_{l_*} \notin H_3$. If $n\in B$ and $\rho \in Y_n$,
then $f(\rho)=F_1(f_{l_n})(\rho)=F_1(f_{l_*})(\rho)$, therefore,
$eq(f,F_1(f_{l_*}))$ is infinite. As $f_{l_*} \notin H_3$, $F_1(f_{l_*}) \in G_4$
(by the definition of $G_4$). Therefore, we've found $g\in G_4$
such that $eq(f,g)$ is infinite and we're done.

Subcase II: $f_{l_*} \in H_3$. For each $n\in B$, $Y_n$ is infinite,
therefore we can find $\rho_n \in Y_n \setminus w_{F_3^*(f_{l_*})}$
(by the definition of $G_2$, $\{\rho \in w_{F_3^*(f_{l_*})} : \rho \cap \eta_f= \eta_f \restriction n\}$
is finite, and as $Y_n \subseteq \{\rho \in 2^{<\omega} : \rho \cap \eta_f= \eta_f \restriction n\}$
is infinite, there is $\rho_n$ as required).

As $f_{l_*} \in H_3$, $f_{l_*}=f_g$ for some $g\in G_2$, and $F_3^*(f_{l_*})=g$,
hence $F_3^*(f_{l_*})(\rho_n)=g(\rho_n)=F_1(f_{l_*})(\rho_n)=f(\rho_n)$
(the equalities follow from the definitions of $F_1$, $F_3^*$ and
$Y_n$, and the assumption that $\rho_n \notin w_{F_3^*(f_{l_*})}$).
Therefore, $eq(F_3^*(f_{l_*}),f)$ is infinite, and by the definition
of $G_4$, $F_3^*(f_{l_*}) \in G_4$ so we're done. $\square$

\textbf{Proof of claim 6: }For $f\in \mathcal{F}_*$, let $\eta_f=\bold{B}(f)$
and let $TV_*(f)$ be the truth value of the statement:

$(*)$ For every $n<\omega$ there exists $\nu \in 2^{<\omega}$ such
that $\eta_f \restriction n \leq \nu$ and $f(\nu) \notin \mathcal{F}_{lg(\nu)}^*$.

Note that $TV_*(f)$ is Borel-computable and so are the truth values
$TV_{2,k,i}(f)$ and $TV_{3,j}(f)$ (to be defined later), therefore,
it suffices to define $F_3^*$ separately for each combination of
truth values.

\textbf{Case I: $TV_*(f)=true$. }In this case, we shall prove that
$f\in H_3$ and define $F_3^*(f)$: 

Let $A_f$ be the set of $n$ for which there is $\nu \in 2^{<\omega}$
such that $\nu \cap \eta_f=\eta_f \restriction n$ and $f(\nu) \notin \mathcal{F}_{lg(\nu)}^*$.
By the assumption, $A_f$ is infinite.

For each $n\in A_f$, let $\nu_{f,n}$ be a sequence for which $(*)$
is true, such that:

1. $lg(\nu_{f,n})$ is minimal.

2. $\nu_{f,n}$ is $<_*-$minimal among the sequence satisfying (1)
(where $<_*$ is the lexicographic ordering).

Let $w_f=\{\nu_{f,n} : n\in A_f\}$ and let $F_3^*(f)=f\restriction w_f \cup F_1(f) \restriction (2^{<\omega} \setminus w_f)$.
It's straightforward to verify that $F_3^*(f) \in G_2$ and that $f_{F_3^*(f)}=f$,
therefore $f\in H_3$.

\textbf{Case II: $TV_*(f)=false$. }We can compute $m(f)=min\{m : $
If $\eta_f \restriction m \leq \nu \in 2^{<\omega}$ then $f(\nu) \in \mathcal{F}_{lg(\nu)}^*\}$.
Let $TV_{2,k,i}(f)$ be the truth value of the following statement:

$(*)_{2,k,i}$ There exist $k$ and $f_0,...,f_{k-1} \in \mathcal{F}_*$
such that for every $\nu \in 2^{<\omega}$, $\eta_f \restriction i \leq \nu \rightarrow f(\nu) \in \{F_1(f_l)(\nu) : l<k\}$.

By compactness, $(*)_{2,k,i}$ holds iff for every finite $u\subseteq \{\nu : \eta_f \restriction i \leq \nu \in 2^{<\omega}\}$
there exist $f_0,...,f_{k-1}$ as above with domain $2^{<lg(u)+1}$
where $lg(u)=max\{lg(\nu) : \nu \in u\}$. Therefore, $TV_{2,k,i}(f)$
is Borel-computable. Note that there is no essential use of the axiom
of choice in the compactness argument, as we can argue in an appropriate
$L[X]$.

Note that if $TV_{2,k,i}(f)=true$ for some $k$ and $i$, then $f\notin H_3$:
Let $f_0,...,f_{k-1}$ be as in $(*)_{2,k,i}$ and suppose towards
contradiction that there exists $g\in G_2$ such that $f=f_g$. Let
$(\nu_n : n\in A)$ list $w_g$, then one of the two possiblities
in 2(m)(2) holds. As $TV_*(f)=false$, the first possibility of 2(m)(2)
fails. Suppose that the second possibility holds. By 2(j)(1), for
every $n\in A \setminus i$ there is $m(n) \in A \setminus i$ such
that $\bold{B}(f) \restriction n \leq \nu_{m(n)}$. As $TV_{2,k,i}(f)=true$,
for every such $n\in A \setminus i$, there exists $l<k$ such that
$f(\nu_{m(n)})=F_1(f_l)(\nu_{m(n)})$. Therefore, for some $l_*<k$,
the set $B:=\{n\in A \setminus i : f(\nu_{m(n)})=F_1(f_{l_*})(\nu_{m(n)})\}$
is infinite. It follows that the elements of $(f(\nu_{m(n)}) : n\in B)=(g(\nu_{m(n)}) : n\in B)$
are pairwise comparable, contradicting the second possibility of 2(m)(2).
It follows that $f\notin H_3$. If $TV_{2,k,i}(f)=true$ for some
$k,i$, we let $F_3^*(f)=f_0$ where $f_0(\nu)$ is defined as the
constant function $0$ whose domain is $2^{<lg(\nu)}$. It's easy
to see that $f_0 \notin G_2$.

From now on, we assume that $TV_*(f)=false$ and $\underset{k,i}{\wedge}TV_{2,k,i}(f)=false$.
As in the proof of claim 8(c) (case III), we shall choose by induction
$A_k$, $h_k$ and $\bar{\eta_k}=(\eta_{k,n} : n_1 \leq n\in A_k)$
(where $n_1=m(f)$ is witnessing the failure of the statement of case
I in the proof of 8(c), and $h_k$ here stand for $f_k$ there) such
that:

a. $A_k \subseteq \omega$ is infinite.

b. $\eta_f \restriction n \leq \eta_{k,n}$ and $\eta_f \restriction (n+1) \nleq \eta_{k,n}$.

c. $(f(\eta_{k,n}) : n_1 \leq n \in A_k)$ is $\subseteq-$increasing.

d. $h_k=\underset{n\in A_k}{\cup}f(\eta_{k,n}) \in \mathcal{F}_*$.

e. $f(\eta_{k,n}) \notin \{F_1(h_l)(\eta_{k,n}) : l<k\}$.

Moreover, the objects will be computed in a Borel way, The only non-trivial
point is the application of Ramsey's theorem in the construction of
$A_k$ from $A_k^1$ (i.e. why can we Borel-compute an infinite homogeneous
set?): Given a function $R: [\omega]^2 \rightarrow \{0,1\}$, we shall
Borel-compute an infinite homogeneous set (we shall write $R(m,k)$
for $R(\{m,k\})$ where $m<k$). Define $\rho_n \in 2^n$ by induction
on $n$ such that:

a. $\rho_n \leq \rho_{n+1}$.

b. For infinitely many $k<\omega$, $R(m,k)=\rho_{n+1}(n)$ for every
$m<n+1$. Let $A_n$ be the set of these $k$'s.

c. $A_{n+1} \subseteq A_n$.

d. $\rho_{n+1}(n)=0$ if possible (i.e. if the above requirements
are satisfied).

The sequence $(\rho_n : n<\omega)$ can be Borel-computed. Now choose
$n_i \in \omega$ by induction such that:

a. $n_i<n_{i+1}$.

b. $n_i$ is the minimal $k \in A_{n_{i-1}}$ such that $\underset{j<i}{\wedge}n_j<k$
and $R(n_j,k)=\rho_{n_j+1}(n_j)$ (this is possible by the choice
of the $\rho_n$s).

So $(n_i : i<\omega)$ is Borel-computable as well. If $i_1<i_2<i_3$
then $R(n_{i_1},n_{i_2})=\rho_{n_{i_1}+1}(n_{i_1})=R(n_{i_1},n_{i_3})$.
Let $i(*) \in \{0,1\}$ be the minimal such that $\{i : R(n_i,n_{i+1})=i(*)\}$
is infinite (this is Borel-computable as well). Finally, the set $\{n_i : R(n_i,n_{i+1})=i(*)\}$
is a Borel-computable infinite homogeneous set. This completes the
argument on the induction.

Let $TV_{3,j}(f)$ be the truth value of the following statement (which
is Borel-computable):

$(*)_{3,j}$ The set $\{h_k \restriction 2^{<j} : k<\omega\}$ is
infinite.

\textbf{Case I: $\underset{j<\omega}{\vee}TV_{3,j}(f)=true$. }In
this case we can Borel compute $F_3^*(f)$ witnessing that $f\in H_3$.
Let $j_f$ be the minimal $j$ such that $TV_{3,j}(f)=true$ (this
is Borel-computable) and let $B=\{k: j_f<k, h_k \restriction 2^{<j_f} \notin \{h_l \restriction 2^{<j_f} : l<k\}\}$,
this set is infinite by our assumption. We choose $n_k \in A_k \setminus j_f$
by induction on $k\in B$ such that $k<k' \rightarrow n_k<n_{k'}$.
Let $w_f=\{\eta_{k,n_k} : k\in B\}$ and let $F_3^*(f)=f \restriction w_f \cup (F_1(f) \restriction 2^{<\omega} \setminus w_f)$.
It's easy to check that $F_3^*(f) \in G_2$ is witnessed by $f$ (hence
$f\in H_3$): If there are $k\neq k' \in B$ such that $F_3^*(f)(\eta_{k,n_k})=F_1(f)(\eta_{k,n_k})$
and $F_3^*(f)(\eta_{k',n_{k'}})=F_1(f)(\eta_{k',n_{k'}})$, then $f(\eta_{k,n_k})=h_k \restriction 2^{<lg(\eta_{k,n_k})}$
and $f(\eta_{k',n_{k'}})=h_{k'} \restriction 2^{<lg(\eta_{k',n_{k'}})}$
are comparable, contradicting the definition of $B$. It's easy to
check that the other requirements in the definition of $G_2$ are
satisfied as well.

\textbf{Case II: $\underset{j<\omega}{\wedge}TV_{3,j}(f)=false$.
}We can Borel-compute a set $B\in [\omega]^{\omega}$ such that $(h_j : k\in B)$
converges to some $h^* \in \mathcal{F}_*$: Let $B_0=B_1=\omega$.
As $TV_{3,2}(f)=false$, there exists $k(2)>2$ such that for infinitely
many $k$, $h_{k(2)} \restriction 2^{<2}=h_k \restriction 2^{<2}$.
Choose $k(2)$ to be the minimal number with the above property and
let $B_2=\{k\in B_1 : k>k(2)$ and $h_{k(2)} \restriction 2^{<2}=h_k \restriction 2^{<2}\}$.
As $TV_{3,3}(f)=false$, there is a minimal $k(3) \in B_2$ such that
$h_{k(3)} \restriction 2^{<3}=h_k \restriction 2^{<3}$ for infinitely
many $k \in B_2$, let $B_3=\{ k\in B_2 : k>k(3)$ and $h_{k(3)} \restriction 2^{<3}=h_k \restriction 2^{<3}\}$.
We now continue the construction by induction, and obtain the set
$B=\{k(2)<k(3)<k(4)<...\}$. Now let $h^*=\underset{n<\omega}{\cup}h_{k(n)} \restriction 2^{<n}$,
it's easy to see that $B$ and $h^*$ are as required. Note that as
$l\neq k \rightarrow h_l \neq h_k$ (by the definition of the $h_k$s),
there is at most one $k$ such that $h_k=h^*$.

We can Borel-compute $(k_i,n_i,m_i)$ by induction on $i$ such that:

1. $k_i \in B$ is increasing with $i$.

2. $m(f)\leq n_i \in A_{k_i}$ is increasing with $i$.

3. $m_i=lg(\eta_{k_i,n_i})$.

4. $f(\eta_{k_i,n_i}) \nsubseteq h^*$.

5. If $j<i$ then $f(\eta_{k_i,n_i}) \restriction 2^{<m_j} \subseteq h^*$.

The induction step: Suppose that we've carried the induction up to
$i$ and let $j=i-1$. Let $m_{i(*)}=max\{m_l : l<i\}$. As $\underset{n\in B}{lim}h_n=h^*$,
for every $n\in B$ large enough (say, $n_* \leq n$ for some $n_*$)
we have $h_n \restriction 2^{<m_{i(*)}} \subseteq h^*$. Let $k_i \in B$
be the first such $n$ above $k_j$ such that, in addition, $h_{k_i} \neq h^*$
(recall that there is at most one $n$ for which $h_n=h^*$). Recall
that $h_{k_i}=\underset{m(f)\leq n\in A_{k_i}}{\cup}f(\eta_{k_i,n})$,
and for $n_1=m(f)\leq n \in A_{k_i}$ large enough, $h_{k_i} \restriction 2^{<n} \nsubseteq h^*$
(otherwise, $h_{k_i}=h^*$, which is a contradiction).

Let $n_i \in A_{k_i} \setminus n_1$ be the first such $n$ above
$n_j$, and let $m_i=lg(\eta_{k_i,n_i})$, so we have carried the
induction successfully. Now let $w_f=\{\eta_{k_i,n_i} : i<\omega\}$
and let $F_3^*(f)=f\restriction w_f \cup (F_1(f) \restriction 2^{<\omega} \setminus w_f)$.
It's easy to check that $F_3^*(f) \in G_2$ as witnessed by $f$,
which belongs to $H_3$. $\square$

\textbf{\large References}{\large \par}

{[}Br{]} Joerg Brendle, Some problems concerning mad families, note,
2013.

{[}KSZ{]} Bart Kastermans, Juris Steprans and Yi Zhang, Analytic and
co-analytic families of almost disjoint functions, J. Symbolic Logic
73 (2008), 1158-1172.

{[}Ma{]} A. R. D. Mathias, Happy families, Ann. Math. Logic \textbf{12
}(1977), no. 1, 59-111. MR 0491197 .

{[}HwSh:1090{]} Haim Horowitz and Saharon Shelah, Can you take Toernquist's
inaccessible away?, arXiv:1605.02419.

{[}To{]} Asger Toernquist, Definability and almost disjoint families,
arXiv:1503.07577.

$\\$

(Haim Horowitz) Einstein Institute of Mathematics

Edmond J. Safra campus,

The Hebrew University of Jerusalem.

Givat Ram, Jerusalem, 91904, Israel.

E-mail address: haim.horowitz@mail.huji.ac.il

$\\$

(Saharon Shelah) Einstein Institute of Mathematics

Edmond J. Safra campus,

The Hebrew University of Jerusalem.

Givat Ram, Jerusalem, 91904, Israel.

Department of Mathematics

Hill Center - Busch Campus,

Rutgers, The State University of New Jersey.

110 Frelinghuysen road, Piscataway, NJ 08854-8019 USA

E-mail address: shelah@math.huji.ac.il
\end{document}